\documentclass[a4paper,10pt,reqno]{amsart}
\usepackage[a4paper,tmargin=40mm,bmargin=35mm,hmargin=30mm]{geometry}

\usepackage[T1]{fontenc}
\usepackage{lmodern}
\usepackage{amsmath}
\usepackage{amssymb}
\usepackage{amsthm}
\usepackage{stmaryrd}
\usepackage{tikz}
\usetikzlibrary{cd}
\usetikzlibrary{decorations.pathmorphing}
\usepackage{booktabs}
\usepackage{multirow}
\usepackage{multicol}
\usepackage{enumitem}
\usepackage{here}
\usepackage{aliascnt}
\usepackage{hyperref}
\usepackage[noabbrev]{cleveref}

\newcommand{\NewTheorem}[3]{
	\newaliascnt{#1}{TheoremEnvironment}
	\newtheorem{#1}[#1]{#2}
	\aliascntresetthe{#1}
	\crefname{#1}{#2}{#3}
	\Crefname{#1}{#2}{#3}
}

\theoremstyle{definition}

\NewTheorem{definition}{Definition}{Definitions}
\NewTheorem{remark}{Remark}{Remarks}
\NewTheorem{axiom}{Axiom}{Axioms}
\NewTheorem{example}{Example}{Examples}
\NewTheorem{observation}{Observation}{Observations}
\NewTheorem{convention}{Convention}{Conventions}
\NewTheorem{notation}{Notation}{Notations}
\NewTheorem{setting}{Setting}{Settings}
\NewTheorem{question}{Question}{Questions}
\NewTheorem{answer}{Answer}{Answers}
\NewTheorem{conjecture}{Conjecture}{Conjectures}
\NewTheorem{problem}{Problem}{Problems}
\NewTheorem{solution}{Solution}{Solutions}
\NewTheorem{goal}{Goal}{Goals}
\NewTheorem{comment}{Comment}{Comments}
\NewTheorem{aim}{Aim}{Aims}
\NewTheorem{caution}{Caution}{Cautions}
\NewTheorem{exercise}{Exercise}{Exercises}

\theoremstyle{plain}
\NewTheorem{proposition}{Proposition}{Propositions}
\NewTheorem{lemma}{Lemma}{Lemmas}
\NewTheorem{theorem}{Theorem}{Theorems}
\NewTheorem{corollary}{Corollary}{Corollaries}

\crefname{enumi}{}{}
\Crefname{enumi}{}{}
\creflabelformat{enumi}{(#2#1#3)}
\crefname{enumii}{}{}
\Crefname{enumii}{}{}
\creflabelformat{enumii}{(#2#1#3)}
\crefname{enumiii}{}{}
\Crefname{enumiii}{}{}
\creflabelformat{enumiii}{(#2#1#3)}

\makeatletter
\renewcommand{\p@enumii}{}
\renewcommand{\p@enumiii}{}
\makeatother

\numberwithin{equation}{section}
\crefname{equation}{}{}
\Crefname{equation}{}{}
\creflabelformat{equation}{(#2#1#3)}

\newcommand{\SwapSymbols}[1]{
	\expandafter\let\expandafter\temporarysymbol\csname #1\endcsname
	\expandafter\let\csname #1\expandafter\endcsname\csname var#1\endcsname
	\expandafter\let\csname var#1\endcsname\temporarysymbol
}

\SwapSymbols{epsilon}
\SwapSymbols{phi}
\SwapSymbols{Gamma}
\SwapSymbols{Delta}
\SwapSymbols{Theta}
\SwapSymbols{Lambda}
\SwapSymbols{Xi}
\SwapSymbols{Pi}
\SwapSymbols{Sigma}
\SwapSymbols{Upsilon}
\SwapSymbols{Phi}
\SwapSymbols{Psi}
\SwapSymbols{Omega}

\newcommand{\bbC}{\mathbb{C}}

\newcommand{\cC}{\mathcal{C}}

\newcommand{\cF}{\mathcal{F}}
\newcommand{\cG}{\mathcal{G}}

\newcommand{\cS}{\mathcal{S}}
\newcommand{\cT}{\mathcal{T}}

\newcommand{\kc}{\mathfrak{c}}

\let\originalleft\left
\let\originalright\right
\renewcommand{\left}{\mathopen{}\mathclose\bgroup\originalleft}
\renewcommand{\right}{\aftergroup\egroup\originalright}

\newcommand{\setwithcondition}[3][]{\mathopen{#1\{}\,#2\mathrel{#1|}#3\,\mathclose{#1\}}}

\DeclareMathOperator{\Hom}{Hom}
\DeclareMathOperator{\End}{End}

\DeclareMathOperator{\Ext}{Ext}

\let\Im\relax
\DeclareMathOperator{\Im}{Im}

\DeclareMathOperator{\Mod}{Mod}

\DeclareMathOperator{\proj}{proj}

\DeclareMathOperator{\add}{add}

\DeclareMathOperator{\Frac}{Frac}

\DeclareMathOperator{\Tr}{Tr}

\DeclareMathOperator{\cHom}{\mathscr{H}\kern -2.75pt\mathit{om}}
\DeclareMathOperator{\cEnd}{\mathscr{E}\kern -1.5pt\mathit{nd}}
\DeclareMathOperator{\cExt}{\mathscr{E}\kern -1.5pt\mathit{x\kern -0.25pt t}}
\DeclareMathOperator{\cTor}{\mathscr{T}\kern -3.25pt\mathit{or}}

\newcommand{\ABfourst}{\mathrm{AB4}^{\ast}}
\newcommand{\ABsix}{\mathrm{AB6}}
\newcommand{\np}{\mathfrak{n}_{+}}
\newcommand{\nm}{\mathfrak{n}_{-}}

\title[A Grothendieck category with a noetherian generator and $\ABfourst$]{A Grothendieck category with a noetherian generator and exact products that is not a module category}

\author{Ryo Kanda}
\address[Ryo Kanda]{Department of Mathematics, Graduate School of Science, Osaka Metropolitan University, 3-3-138, Sugimoto, Sumiyoshi, Osaka, 558-8585, Japan}
\email{ryo.kanda.math@gmail.com}


\begin{document}

\begin{abstract}
We construct a Grothendieck category that has a noetherian generator, satisfies $\ABfourst$, and is not equivalent to a module category. This gives a negative answer to Djament's problem. The category is obtained as a Gabriel quotient of a module category over the endomorphism ring appearing in the work of Herbera, P\v{r}\'{\i}hoda, and Wiegand. The proof relies on the trace criterion established by Martini, Parra, Saor\'{\i}n, and Virili.
\end{abstract}

\maketitle

\section{Introduction}
\label{58310472}

An abelian category with direct products is said to satisfy $\ABfourst$ if its direct product functors are exact. Martini--Parra--Saor\'{\i}n--Virili recorded the following problem, posed by Aur\'elien Djament:

\begin{problem}[{\cite[Problem~5.13]{arXiv:2508.00670v2}}]\label{94627138}
Let $\cG$ be a locally noetherian, locally coherent, or locally finitely presented Grothendieck category satisfying $\ABfourst$. Does $\cG$ have a generating set consisting of finitely generated projective objects?
\end{problem}

Independently, the author recorded in \cite[Remark~7.15]{MR4342508} that he was not aware of any example of a Grothendieck category with a noetherian generator and exact products that is not equivalent to a module category.

Martini--Parra--Saor\'{\i}n--Virili \cite[Theorem~5.14 and Proposition~5.15; see also Remark~5.16(4)]{arXiv:2508.00670v2} gave negative answers in the locally coherent and locally finitely presented cases. More precisely, they constructed a locally coherent Grothendieck category with exact products and enough flat objects, but without enough projectives. Since a locally coherent Grothendieck category is locally finitely presented, this also gives a negative answer in the locally finitely presented case. They note that the locally noetherian case remains open; see \cite[Introduction]{arXiv:2508.00670v2}. The purpose of this paper is to give a negative answer in this remaining case. In fact, we construct such an example with a noetherian generator. This is precisely the kind of example sought in \cite[Remark~7.15]{MR4342508}.

The construction starts from the local ring $R$ of a nodal curve singularity and its normalization $\overline{R}$. We consider the endomorphism ring $T=\End_{R}(R\oplus \overline{R})$, which appears in the work of Herbera--P\v{r}\'{\i}hoda--Wiegand \cite[Section~8]{MR4915561} on pure projective modules over commutative noetherian rings. The normalization separates the node into two branches. Choosing one branch gives an idempotent two-sided ideal $I\subseteq T$. The Gabriel quotient associated with this ideal is the desired Grothendieck category.

The main point is that this quotient has exact products and a noetherian generator, while the ideal $I$ is not the trace of any set of finitely generated projective right $T$-modules. By the trace criterion of Martini--Parra--Saor\'{\i}n--Virili \cite[Theorem~5.7]{arXiv:2508.00670v2} (see \cref{42681573} below), the quotient therefore does not admit a generating set consisting of finitely generated projective objects.

\begin{theorem}\label{81736405}
There exists a Grothendieck category satisfying $\ABfourst$ and having a noetherian generator that does not admit a generating set consisting of finitely generated projective objects.
\end{theorem}

Such a category is not equivalent to the module category over a small preadditive category (in particular, over a ring), since module categories over small preadditive categories admit a generating set consisting of finitely generated projective objects (see \cite[Section~2.2]{arXiv:2508.00670v2}). A Grothendieck category with a noetherian generator is locally noetherian, so this gives a negative answer to the locally noetherian case of Djament's problem.

Since the example we construct has a projective generator that is not finitely generated (see \cref{47160832}), the following problem remains open:

\begin{problem}
\label{93282996}
Is there a locally noetherian Grothendieck category satisfying $\ABfourst$ that does not admit a projective generator?
\end{problem}

\subsection*{Tool and Computational Resource Disclosure}

The author used OpenAI's ChatGPT, with models from the GPT-5.5 and GPT-5.6 families, as a general-purpose tool for research assistance and manuscript preparation, including exploratory assistance and suggestions concerning exposition and language revision. In particular, the tool suggested examining the example in \cite[Section~8]{MR4915561} as a possible source of a negative answer to \cite[Problem~5.13]{arXiv:2508.00670v2}. The tool also suggested a possible proof strategy. The author independently checked this strategy and supplied the proof presented in this paper. The author is solely responsible for the content of this paper.

\subsection*{Acknowledgments}

The author was supported by JSPS KAKENHI Grant Number JP24K06693, the MEXT Promotion of Distinctive Joint Research Center Program JPMXP0723833165, and the Osaka Metropolitan University Strategic Research Promotion Project (Development of International Research Hubs).

\section{The trace criterion of Martini--Parra--Saor\'{\i}n--Virili}
\label{69412037}

The key result used in this paper is the trace criterion of Martini--Parra--Saor\'{\i}n--Virili \cite[Theorem~5.7]{arXiv:2508.00670v2}. Although the original statement is formulated for a small preadditive category, we state only its specialization to a ring, that is, to the one-object case.

\subsection{Idempotent ideals and associated quotient categories}
\label{15278364}

Let $T$ be a ring and let $I\subseteq T$ be an idempotent two-sided ideal. Denote by $\Mod T$ the category of right $T$-modules, and set $\cT_{I}=\{M\in\Mod T\mid MI=0\}$. This category appears in the TTF triple $(\cC_{I},\cT_{I},\cF_{I})$ associated with $I$; see \cite[Proposition~4.1]{arXiv:2508.00670v2}. The Giraud subcategory
\begin{equation*}
\cG_{I}:=\setwithcondition{M\in\Mod T}{\textnormal{$\Hom_{T}(X,M)=0=\Ext^{1}_{T}(X,M)$ for all $X\in\cT_{I}$}}
\end{equation*}
is equivalent to the Gabriel quotient $(\Mod T)/\cT_{I}$ (see \cite[Section~1.2]{arXiv:2508.00670v2}). The category $\cG_{I}$ is a Grothendieck category satisfying $\ABsix$ and $\ABfourst$ (see \cite[Proposition~5.2]{arXiv:2508.00670v2}). If $T$ is right noetherian, then $T_{T}$ is a noetherian generator of $\Mod T$. Since the quotient functor sends noetherian objects to noetherian objects and sends generating sets to generating sets (see, for example, \cite[Proposition~2.6]{MR3922832}), the image of $T_{T}$ is a noetherian generator of $(\Mod T)/\cT_I\cong\cG_I$.

Let $S$ be a right $T$-module. The \emph{trace} of $S$ is the two-sided ideal
\begin{equation*}
\Tr_{T}(S):=\sum_{f\colon S\to T}\Im f
\end{equation*}
where $f$ runs over all right $T$-module homomorphisms $S\to T$. For a class $\cS$ of right $T$-modules, the \emph{trace} of $\cS$ is
\begin{equation*}
\Tr_{T}(\cS):=\sum_{S\in\cS}\Tr_{T}(S).
\end{equation*}
Note that $\Tr_{T}(S_{1}\oplus S_{2})=\Tr_{T}(S_{1})+\Tr_{T}(S_{2})$ for $S_{1},S_{2}\in\Mod T$. For an idempotent $e\in T$, we have $\Tr_{T}(eT)=TeT$, because every right $T$-module homomorphism $eT\to T$ is given by left multiplication by an element of $Te$.

\begin{theorem}[Martini--Parra--Saor\'{\i}n--Virili; {\cite[Theorem~5.7]{arXiv:2508.00670v2}}]\label{42681573}
Let $T$ be a ring, let $I\subseteq T$ be an idempotent two-sided ideal, and let $\cG_{I}\subseteq\Mod T$ be the Giraud subcategory associated with $I$. Then:
\begin{enumerate}
\item The category $\cG_{I}$ has a projective generator if and only if $I$ is the trace of a projective right $T$-module.
\item The category $\cG_{I}$ has a generating set of finitely generated projective objects if and only if $I$ is the trace of a set of finitely generated projective right $T$-modules.
\end{enumerate}
\end{theorem}

Thus, to construct a Grothendieck category that admits a noetherian generator, satisfies $\ABfourst$, but does not admit a generating set consisting of finitely generated projective objects, it suffices to construct a right noetherian ring $T$ and an idempotent ideal $I\subseteq T$ such that $I$ is not the trace of any set of finitely generated projective right $T$-modules.

According to \cite[Theorems~6.1 and~6.4]{arXiv:2508.00670v2}, such a ring $T$ should not be sought among commutative or semiregular rings.

\section{Construction of the example}
\label{36821497}

The ring $T$ used below is taken from the example of Herbera--P\v{r}\'{\i}hoda--Wiegand \cite[Section~8]{MR4915561}, which appears in their study of pure projective modules. We construct an idempotent ideal $I\subseteq T$ and prove that the associated Giraud subcategory gives the desired example.

Consider the local ring of a nodal curve singularity
\begin{equation*}
R=
\biggl(
\dfrac{\bbC[x,y]}{(x^{2}-y^{3}-y^{2})}
\biggr)_{(x,y)}.
\end{equation*}
Put $u=x/y\in\Frac R$. The relation $x^{2}=y^{2}(y+1)$ gives $u^{2}=y+1$, so $y=u^{2}-1$ and $x=u(u^{2}-1)$. The normalization of $R$ is $\overline{R}=R+uR=R[u]\subseteq\Frac R=\bbC(u)$, which is a semilocalization of $\bbC[u]$ at the two maximal ideals $(u-1)$ and $(u+1)$. Denote the corresponding maximal ideals of $\overline{R}$ by $\np=(u-1)\overline{R}$ and $\nm=(u+1)\overline{R}$.

Let $\kc=(R:\overline{R})\subseteq R$ be the conductor ideal. Since $\overline{R}=R+uR$, we have $\kc=\{r\in R\mid ru\in R\}$ and $y\overline{R}=yR+yuR=(x,y)R$. Moreover, $yu=x\in R$ and $xu=y(y+1)\in R$, so $(x,y)R\subseteq\kc$. On the other hand, $u\notin R$, so $1\notin\kc$. Thus $\kc$ is a proper ideal of the local ring $R$, and hence $\kc\subseteq(x,y)R$. Therefore
\begin{equation*}
\kc=(x,y)R=y\overline{R}=(u^{2}-1)\overline{R}=\np\nm
\end{equation*}
and $\kc\overline{R}=\kc$. Set
\begin{equation*}
T=\End_{R}(R\oplus\overline{R})=
\begin{pmatrix}
R & \kc\\
\overline{R} & \overline{R}
\end{pmatrix}.
\end{equation*}
The ring $T$ is naturally an $R$-algebra, with $R$ contained in the center of $T$. Since $\overline{R}$ is finitely generated as an $R$-module, the ring $T$ is also finitely generated as an $R$-module. Since $R$ is noetherian, $T$ is a left and right noetherian ring.

Let $e=\begin{pmatrix}1 & 0 \\ 0 & 0\end{pmatrix}$ and $f=\begin{pmatrix}0 & 0 \\ 0 & 1\end{pmatrix}$ be the idempotents in $T$. The trace of $eT$ is
\begin{equation*}
\Tr_{T}(eT)=TeT=
\begin{pmatrix}
R & \kc\\
\overline{R} & \kc
\end{pmatrix}.
\end{equation*}
The quotient $T/TeT$ is identified with $\overline{R}/\kc$. Define $I\subseteq T$ to be the inverse image of the ideal $\np/\kc\subseteq\overline{R}/\kc$. Explicitly,
\begin{equation*}
I=
\begin{pmatrix}
R & \kc\\
\overline{R} & \np
\end{pmatrix}.
\end{equation*}
Then
\begin{equation*}
I^{2}=
\begin{pmatrix}
R & \kc\\
\overline{R} & \kc+\np^{2}
\end{pmatrix}.
\end{equation*}
Since $\kc=\np\nm$ and $\np+\nm=\overline{R}$, we have $\kc+\np^{2}=\np\nm+\np^{2}=\np$. Therefore $I^{2}=I$.

It remains to show that $I$ is not the trace of any set of finitely generated projective right $T$-modules. As in \cite[Section~8]{MR4915561}, all finitely generated torsionfree $R$-modules are precisely the finite direct sums of copies of $R$ and $\overline{R}$. Thus they form the additive closure $\add(R\oplus\overline{R})$. The equivalence $\Hom_{R}(R\oplus\overline{R},-)\colon\add(R\oplus\overline{R})\to\proj T$ implies that all finitely generated projective right $T$-modules are precisely the finite direct sums of copies of $eT$ and $fT$. We have already computed the trace of $eT$, and the trace of $fT$ is
\begin{equation*}
\Tr_{T}(fT)=TfT=
\begin{pmatrix}
\kc & \kc\\
\overline{R} & \overline{R}
\end{pmatrix}.
\end{equation*}

Therefore the trace of a finitely generated projective right $T$-module is one of the four ideals $0$, $TeT$, $TfT$, and $TeT+TfT=T$. Since the trace of a set of finitely generated projective modules is the sum of the traces of its members, the trace of any such set is again one of the four ideals above. Hence $I$ is not the trace of any set of finitely generated projective right $T$-modules.

As explained in \cref{69412037}, the Giraud subcategory $\cG_{I}\subseteq\Mod T$ is a Grothendieck category that has a noetherian generator and satisfies $\ABfourst$, but does not admit a generating set consisting of finitely generated projective objects by \cref{42681573}. Thus \cref{81736405} is proved.

\begin{remark}\label{47160832}
The Grothendieck category $\cG_{I}$ constructed here has a projective generator that is not finitely generated. Indeed, $I$ is finitely generated as a left ideal of $T$, since $T$ is left noetherian. Since $I$ is also idempotent, Whitehead's theorem gives a countably generated projective right $T$-module $P$ such that $\Tr_{T}(P)=I$; see \cite[Corollary~2.7]{MR588450} and also \cite[Corollary~2.7]{MR3232793}. Applying \cref{42681573}, the associated Giraud subcategory $\cG_{I}$ has a projective generator.
\end{remark}

\bibliography{references}
\bibliographystyle{customamsalpha}

\end{document}